\documentclass[graybox]{svmult}

\usepackage{mathptmx}       % selects Times Roman as basic font
\usepackage{helvet}         % selects Helvetica as sans-serif font
\usepackage{courier}        % selects Courier as typewriter font
\usepackage{type1cm}        % activate if the above 3 fonts are
\usepackage{makeidx}         % allows index generation
\usepackage{graphicx}        % standard LaTeX graphics tool
\usepackage{multicol}        % used for the two-column index
\usepackage[bottom]{footmisc}% places footnotes at page bottom

\usepackage{amssymb}
\usepackage{psfrag}

\usepackage{natbib}

\newcommand\R{\mathbb{R}}
\newcommand\N{\mathbb{N}}
\newcommand\wP{\widetilde P}
\newcommand\wpsi{\widetilde  \psi}
\newcommand{\two}[2]{\begin{array}{c}#1\cr#2\end{array}}
\makeindex             % used for the subject index
\newtheorem{algorithm}[lemma]{Algorithm}

\begin{document}

\title*{A Reduced Basis Method for the Simulation of American Options}
\author{Bernard Haasdonk, Julien Salomon and Barbara Wohlmuth}
\institute{Bernard Haasdonk \at IANS, Universit\"at Stuttgart, Germany,  \email{haasdonk@mathematik.uni-stuttgart.de}
\and Julien Salomon \at CEREMADE, Universit\'e Paris-Dauphine, \email{salomon@ceremade.dauphine.fr}
\and Barbara Wohlmuth \at  M2 - Zentrum Mathematik, Technische Universit\"at M\"unchen, \email{wohlmuth@ma.tum.de}}
\maketitle

%\abstract*{}

\abstract{We present a reduced basis method for the
  simulation of American option pricing. To tackle this model
  numerically, we formulate the problem in terms of a time dependent
  variational inequality. Characteristic ingredients are a
  POD-greedy and an angle-greedy procedure for the construction of 
  the primal and dual reduced spaces.
  Numerical examples are provided, illustrating the approximation quality and convergence of our approach.}

\section{Introduction}\label{sec:1}

We consider the problem of American option pricing and refer to~\citep{achdou} and the references therein for an introduction into
computational methods for option pricing.
While European options can be modelled by a parabolic partial differential
equation,
American options result in additional inequality constraints.
We refer to~\citep{wohlmuth}  for a possible numerical
treatment by primal-dual finite elements and to~\citep{glowinski,geiger} for an abstract framework on the theory of constrained variational
problems. 
We are interested in providing a fast numerical algorithm to solve
accurately the variational
 inequality system of an American put option for a large variety of
different parameter values
such as interest rate, dividend, strike prize and volatility.
Reduced basis (RB) methods are an appropriate means
for standard para\-metrized parabolic partial differential equations,
cf. \citep{haasdonk08,rozza,veroy,buffa}
and the references therein. These are based on
low-dimensional approximation spaces, that are constructed
by greedy procedures. Convergence behavior of these procedures
are known in some cases
\citep{buffa, haasdonk11}.
The computational advantage of RB-methods over standard discretization 
methods 
is obtained by its possible offline/online decomposition: First, a typically 
expensive offline-phase involving the computation of the reduced spaces is 
performed. This phase only needs to be precomputed once. Then, the online phase allows an extremely fast computation of the RB
solutions for many new parameters 
as only low dimensional systems need to be solved. Recently, we adopted the RB methodology to
constrained stationary elliptic problems \citep{haasdonksalomonwohlmuth}, which
we extend here to the instationary case.
We refer to the recent contribution~\citep{pironneau} for a
tailored RB approach in
option pricing. In contrast to our setting no inequality constraints are
taken into account. The main challenge is the construction of a suitable low dimensional
approximation of the dual cone required for
the approximation of the constraints. In this contribution, we introduce a
new greedy
 strategy based on an angle criteria and show numerical results.

\section{American Option Model}\label{sec:2}

An American option is a contract which permits its owner to receive a
certain payoff
$\psi(S, \tau ) \geq 0$
at any time $\tau$ between $0$ and $T>0$.
The variable $T$
indicates the maturity.
Introducing the backward time variable $t:= T - \tau$, we can use, e.g.,~\citep{achdou}
the following non linear model
\begin{eqnarray*}
\partial_t P -\frac 12 \sigma^2 s^2\partial^2_{ss}P-(r-q)s\partial_s P
+r P\geq 0,\quad P-\psi &\geq& 0, \\%\label{eq:original2}\\
\left( \partial_t P -\frac 12 \sigma^2 s^2\partial^2_{ss}P-(r-q)s\partial_s P
+r P\right)\cdot \left(P-\psi \right)&=&0,%\label{eq:original3}
\end{eqnarray*}
where $P=P(s,t)$ is the
price of an American put, with $s\in \R_+$ the asset's
value, $\sigma$ is the volatility, $r$ is the interest rate,
$q$ is the dividend payment and $\psi=\psi(s,t)$ is the payoff function. 
The boundary and initial conditions are as follows:
%\begin{equation}
$P(s,0)=\psi(s),$  $P(0,t)=K,$ $\lim_{s\rightarrow+\infty}P(s,t)=0,$
%\nonumber %\end{equation}
where $K>0$ is a fixed strike price that satisfies $K=\psi(0,0)$. In what
follows, we focus on the case $\psi(s,t )=(K-s)_+$ with
$(\cdot)_+=\max(0,\cdot)$, but our method applies as well to other types of
payoff functions. For the implementation, we restrict the values of $s$ to a
bounded interval $\Omega:=(0,s_f)$, where $s_f$ is large enough to make the
assumption $P(s_f,t)=0$ realistic. Let us also set $\wP = P-P_0$, with
initial data $P_0(s,t)=K(1- s/s_f)$, so that $\wP$ satisfies homogeneous Dirichlet conditions. 
Our aim is now to reformulate the last system in a weak form, where
our reduced basis method applies. In this view, we introduce the following functional spaces:
$$ V:= \left\{ v\in L^2(\Omega) | s\partial_s v\in L^2(\Omega),  v_{|\partial \Omega}=0\right\}, \quad W:=V'.$$
The scalar product $\langle \cdot,\cdot \rangle_V$  associated with $V$
is defined by $\langle u,v \rangle_V:=\langle s\partial_s u,s\partial_s v
\rangle_{L^2(\Omega)}+\langle u,v \rangle_{L^2(\Omega)},$
where $\langle \cdot,\cdot \rangle_{L^2(\Omega)}$ is
the usual scalar product on $L^2(\Omega)$. The operators are specified as follows:
$$ a(u,v;\mu)= \frac 12\sigma^2\langle  \partial_{s}u,
\partial_{s}(s^2v)\rangle_{L^2(\Omega)}+\langle -(r-q)s\partial_su +ru,v\rangle_{L^2(\Omega)},$$
$$f(v;\mu)=\langle F , v\rangle_{L^2(\Omega)}, \ g(\eta ;\mu)=\langle \wpsi ,\eta \rangle_W, $$
with  $F:= -\left( \partial_t P_0 -\frac 12 \sigma^2 s^2\partial^2_{ss}P_0-(r-q)s\partial_s P_0 +r P_0\right)$, i.e.
$F=K\left(\frac s {s_f} q-r\right)$ and $\wpsi:=\psi-P_0$. For $\eta\in W=V'$, we also define $b(\eta,v)= \eta(v)$.
We can now recast our problem  in the following weak
form, para\-me\-trized by  $\mu=(K,r,q,\sigma)\in {\cal P}\subset
\R^4$. We now introduce $u$ as a weak representant of 
the solution $\widetilde P$, as 
this is the standard notation in reduced basis literature: 
\begin{eqnarray}
\langle\partial_t u,v\rangle_{L^2(\Omega)}+a(u,v;\mu)-b(\lambda,v) & =& f(v;\mu),
\qquad v\in V\label{eq:VIgene1} \\
b(\eta-\lambda,u) & \geq& g(\eta-\lambda;\mu), \qquad \eta\in M,\label{eq:VIgene2}
\end{eqnarray}
where $M\subset W$ is a closed convex cone.
Various methods can be considered to solve numerically
Equations~(\ref{eq:VIgene1}--\ref{eq:VIgene2}). In what follows, we use a
$\theta$-scheme for the time discretization. Given $\mu\in\cal
P$,  $L\in\N$ and
$\Delta t:=T/L$, this method
corresponds to the following iteration.

Given $0<n\leq L-1$ and $u^n\in V$,
  find $u^{n+1}\in V$ and $\lambda^{n+1}\in M$ that satisfy $\forall
  v\in V, \forall \eta \in M$,
\begin{eqnarray}
\left\langle\!\!\frac{u^{n+1}-u^n}{\Delta t},v\right\rangle_{L^2(\Omega)}
\!\!\!\!\!\!\!\!\!\!\!\!+a(\theta
  u^{n+1}+(1-\theta)u^n,v;\mu) - b(\lambda^{n+1},v) & = & f(v;\mu),
\label{eq:VIdisc1}  \\ 
b(\eta-\lambda^{n+1},u^{n+1}) & \geq & g(\eta-\lambda^{n+1};\mu).\quad
 \label{eq:VIdisc2} 
\end{eqnarray}
This recursive definition is initialized with $u^0:=\widetilde \psi$. Note that
in this scheme, the definition of $\lambda^n$ is not recursive. 

\section{Reduced Basis Method}\label{sec:3}
Standard finite element approaches do not exploit the structure of the
solution and for a given parameter value, a high dimensional system has to be solved.
In what follows, we introduce a
specific Galerkin approximation of the solution, based on the reduced
basis method and present algorithms to compute the corresponding bases. 
%\subsection{Reduced basis formulation}
The principle of
the reduced basis method consists in computing parametric  
solutions in low dimensional subspaces of $V$ and
$W$ that are generated with particular solutions of our
problem. Let us explain in more detail the corresponding formulation.
For $N\in\N$, consider a finite subset ${\cal
P}_N:=\left\{\mu_1,\dots,\mu_N \right\}\subset {\cal P}$ with
$\mu_i\neq\mu_j,\ \forall i\neq j$. The reduced spaces
$V_N$ and $W_N$ are defined by $V_N:={\rm
  span}\left\{\psi_1,\dots,\psi_{N_V}\right\}$ and $W_N:={\rm
  span}\left\{\xi_1,\dots,\xi_{N_W}\right\}$ 
%\begin{eqnarray*}
%V_N&:={\rm span}\left\{\psi_1,\dots,\psi_{N_V}\right\}
%\subset {\rm span}\left\{       u^n(\mu_i),\ i=1,\dots,N,\ n=0,\dots,L\right\},\\
%W_N&:={\rm span}\left\{\xi_1,\dots,\xi_{N_W}\right\}
%\subset {\rm span}\left\{ \lambda^n(\mu_i),\ i=1,\dots,N,\ n=0,\dots,L\right\},
%\end{eqnarray*}
%{\bf To be checked again...}
where $\psi_i$ and $\xi_i$ are defined from
the large set of snapshot solutions $u^n(\mu_i)$ and $\lambda^n(\mu_i)$, $i= 1, \ldots, N
$, $n = 0, \ldots, L$. Here $u^n(\mu_i) $ and $\lambda^n (\mu_i)$
denote the solution of Equations~(\ref{eq:VIdisc1}--\ref{eq:VIdisc2}) at the time $t_n := n \Delta t$ for
the parameter value  $\mu = mu_i$.
The
functions $\psi_j$ and $\xi_j$ are suitably selected elements spanning
$V_N$ and $W_N$ with  $N_V,N_W\leq N(L+1)$ preferably small.
Both
families $\Psi_N=(\psi_j)_{j=1,\dots,N_V}$ and $\Xi_N=(\xi_j)_{j=1,\dots,N_W}$ are supposed to be composed of
linearly independent functions, hence are so called reduced bases.
%Because of possible linear dependency of $\xi_j$ The family $\Xi_N=(\xi_j)_{j=1,\dots,N_W}\subset M$ is not necessarily
%of dimension $N_W$ but one has dim $W_N\leq N_W$. 
Numerical algorithms to build these two sets will be presented in Section~\ref{sec:RBcomp}. 
We define the reduced cone $M_N\subset M$ as
$$M_N=\left\{ \sum_{j=1}^{N_W} \alpha_j\xi_j,\ \alpha_j\geq
0\right\}.$$
In this setting, the reduced problem reads:

 Given $\mu\in{\cal
P}$, $0\leq  n\leq L-1$, $u_N^n\in V_N$,
  find $u_N^{n+1}\in V_N$ and $\lambda_N^{n+1}\in M_N$ that satisfy $\forall
  v_N\in V_N, \forall \eta_N \in M_N$,
\begin{eqnarray}
\left\langle\frac{u_N^{n+1}-u_N^n}{\Delta t},v_N\right\rangle_{L^2(\Omega)} \!\!\!\!\!\!\!\!\!\!\!\!+a(\theta
  u_N^{n+1}+(1-\theta)u_N^n,v_N;\mu) - b(\lambda_N^{n+1},v_N) = f(v_N;\mu), \label{eq:VIRB1}  \\ 
b(\eta_N-\lambda_N^{n+1},u_N^{n+1}) \geq g(\eta_N-\lambda_N^{n+1};\mu), \label{eq:VIRB2} 
\end{eqnarray}
where the initial value $u_N^0$ is chosen as the orthogonal projection
of $u_0$ on $V_N$, i.e. 
$$ \langle u_N^0-u_0,v_N \rangle_V=0,\ \forall v_N\in V_N.$$

\section{Reduced Basis Construction}
\label{sec:RBcomp}
In this section, we present two methods to extract a basis $\Psi_N\subset V$ and $\Xi_N \subset M$ from the snapshots. 
Both are greedy procedures based on a finite training set
${\cal P}_{train}\subset {\cal P}$ 
small enough such that it can be scanned quickly.
Given an arbitrary 
integer $N_W$, the dual reduced basis
$\Xi_N=(\xi_j)_{j=1,\dots,N_W}$ is built iteratively according to the
following algorithm. The goal of the approach is to obtain a reduced cone 
$M_N\subset M$ capturing as much ``volume'' as possible. 
\begin{algorithm}(Angle-greedy algorithm)\label{alg:ang-greedy}
Given $N_W$, ${\cal P}_{train}\subset {\cal P}$, choose arbitrarily $0\leq n_1 \leq L$ and $\mu_1\in {\cal P}_{train}$
and do
\begin{enumerate}
\item set $\Xi^1_N=\left\{  \frac{ \lambda^{n_1}(\mu_1) } { \|\lambda^{n_1}(\mu_1)\|_W } \right\}$, $W_N^1:={\rm span}(\Xi_N^1),$
\item for $k=1,\ldots,N_W-1$, do 
\begin{enumerate}
\item find $(n_{k+1},\mu_{k+1}):={\rm argmax}_{n=0,\ldots,L , \ \mu\in {\cal
    P}_{train}} \left( \measuredangle \left( \lambda^n(\mu) , W_N^k\right)\right),$
\item set $\xi_{k+1}:=\frac{\lambda^{n_{k+1}}(\mu_{{k+1}})}{\|\lambda^{n_{k+1}}(\mu_{{k+1}})\|_W}$,
\item define $\Xi^{k+1}_N=\Xi^k_N\cup \{\xi_{k+1}\}$, $W_N^{k+1}:={\rm
span}(\Xi_N^{k+1})$,
\end{enumerate}
\item define $\Xi_N:=\Xi_N^{N_W}$, $W_N:={\rm span}(\Xi_N)$.
\end{enumerate}
\end{algorithm}
Here we have used the notation $\measuredangle (v,S)$ to denote the angle
between a vector $v$ and a linear space $S\subset W$, which is simply 
obtained via the orthogonal projection $\Pi_{S}$ from $W$ on $S$ by
$$
 \measuredangle (v,S) = \arccos \frac{ ||\Pi_{S}v||_W}{ ||v||_W}, \quad v\in W.
$$
We apply the POD-greedy
algorithm ~\citep{haasdonk08}
to design the primal reduced basis $\Psi_N$. 
This procedure is standard in RB-methods 
for evolution problems. 
In RB-methods, frequently {\em weak} greedy procedures are used, 
which make beneficial use of rapidly computable error estimators and 
allow to handle large sets ${\cal P}_{train}$ \citep{buffa}. 
However, as our analysis does 
not yet provide a-posteriori error estimators, we use the true 
projection errors as error indicators. This corresponds to 
the so called {\em strong} greedy procedure \citep{buffa,haasdonk11}.
%It admittedly only allows to use smaller training sets
%${\cal P}_{train}$, but it decouples the space construction 
%from the reduced simulation scheme, and hence focusses
%It has the nice effect of 
%The
%approach consists in enriching the basis by adding
%iteratively a vector representing the trajectory associated to the
%worst approximation. 
%Let us give details about the
%rigorous definitions of this vector and the corresponding trajectory.
\begin{algorithm}(POD-greedy algorithm)\label{alg:pod-greedy}
Given $\widetilde{N}_V>0$, ${\cal P}_{train}\subset {\cal P}$, choose
arbitrarily $\mu_1\in {\cal P}_{train}$,
\begin{enumerate}
\item set $\widetilde\Psi_N^1=\left\{\frac{u^0(\mu_1)}{\|u^0(\mu_1)\|_V}\right\}$,
  $\widetilde V^1_N:={\rm
  span}(\widetilde\Psi_N^1)$,
\item for $k=1,\ldots,\widetilde N_V-1$, do
\begin{enumerate}
\item  define $ \mu_{k+1}:={\rm argmax}_{\mu\in{\cal P}_{train}}{\left(
  \sum_{n=0}^L\|u^n(\mu)-\Pi_{\widetilde V^k_N}(u^n(\mu))\|^2_V   \right)}$,
\item  define $\widetilde\psi_{k+1}:=POD_1\left(u^n(\mu_{k+1})-\Pi_{\widetilde V^k_N}(u^n(\mu_{k+1}))\right)_{n=0,\dots,L},$
\item define $\widetilde\Psi_N^{k+1}:=\widetilde\Psi_N^{k}\cup\left\{\widetilde\psi_{k+1}\right\}$,
\end{enumerate}
\item define $\widetilde\Psi_N:=\widetilde\Psi_N^{ \widetilde N_V}$, $\widetilde V_N:={\rm
  span} \widetilde\Psi_N$.
%define $\widetilde\Psi_N:=\Psi_N^{\widetilde N_V}$ and $V_N^{k+1}:=V_N^k \oplus {\rm span}(\Psi_N^{k+1})$.
\end{enumerate}
\end{algorithm}
Here, we have denoted by $\Pi_{\widetilde V^k_N}$ the orthogonal
projection on $\widetilde V^k_N$ with respect to $\langle \cdot,\cdot\rangle_V$, 
and by $POD_1$ the routine that extracts from a
family of vectors the first Proper Orthogonal
Decomposition (POD) mode that can be obtained via 
the best approximation property
$$
POD_1\left(v^n\right)_{n=0,\dots,L} :=
\arg\min_{||z||_V=1}
  \sum_{n=0}^L|| v^n - \left\langle v^n, z\right\rangle_V z ||_V^2. 
$$
In this definition $V$ is spanned by $v^n$, $n=0,\ldots,L$. A convergence analysis of the POD-greedy
procedure is provided in~\citep{haasdonk11}.
Note that Algorithm~\ref{alg:pod-greedy} always returns an orthonormal basis. This is
even the case if a parameter value
$\mu \in {\cal P}_{train}$ is selected more than once. We point out that our System~(\ref{eq:VIRB1}--\ref{eq:VIRB2}) has a saddle point structure. Thus
taking ${\rm span} \Psi_N$ as
reduced basis for the primal variable
might result in an ill posed problem. To guarantee the inf-sup stability of
our approach, we
follow an idea  introduced in~\citep{rozza} for the Stokes problem, see 
also~\citep{haasdonksalomonwohlmuth}
for variational inequalities. It  
consists in the enrichment $\Psi_N:=\widetilde\Psi_N^{\widetilde N_V}\cup
\left(B\xi_i\right)_{i=1,\ldots,N_W} $, where $B\xi_i$ is the
solution of $b(\xi_i,v)=\langle B\xi_i,v\rangle_V$, for $v\in V$.
We conclude with the final reduced 
space $V_N:={\rm span \Psi_N}$ of dimension $N_V:=$dim $V_N$. 
By construction we have $\widetilde N_V \leq N_V \leq \widetilde N_V + N_W$.

\section{Numerical Results}\label{sec:4}

\begin{figure}
\centering

\psfrag{1.6}[r][m]{\scriptsize 1.6}
\psfrag{0.35}[r][m]{\scriptsize 0.35}
\psfrag{-0.25}[r][m]{\scriptsize -0.25}
\psfrag{0}[r][b]{\scriptsize 0}
\psfrag{150}[c][b]{\scriptsize 150}
\psfrag{300}[c][b]{\scriptsize 300}
\psfrag{40}[r][m]{\scriptsize 40}
\psfrag{80}[r][m]{\scriptsize 80}
\psfrag{20}[c][m]{}
\psfrag{60}[c][m]{}
\psfrag{50}[c][m]{}
\psfrag{100}[c][m]{}
\psfrag{200}[c][m]{}
\psfrag{250}[c][m]{}
\psfrag{-0.3}[c][m]{}
\psfrag{-0.2}[c][m]{}
\psfrag{-0.15}[c][m]{}
\psfrag{-0.1}[c][m]{}
\psfrag{-0.05}[c][m]{}
\psfrag{0.05}[c][m]{}
\psfrag{0.1}[c][m]{}
\psfrag{0.2}[c][m]{}
\psfrag{0.4}[c][m]{}
\psfrag{0.6}[c][m]{}
\psfrag{0.8}[c][m]{}
\psfrag{1}[c][m]{}
\psfrag{1.2}[c][m]{}
\psfrag{1.4}[c][m]{}
\psfrag{0.15}[c][m]{}
\psfrag{0.25}[c][m]{}
\psfrag{0.3}[c][m]{}

\psfrag{p}[c][t]{\footnotesize $\psi_j(s)$}
\psfrag{q}[c][b]{\footnotesize% $s$
$\begin{array}{c} \phantom{1}\\s\end{array}$
}
\includegraphics[width=.26\linewidth]{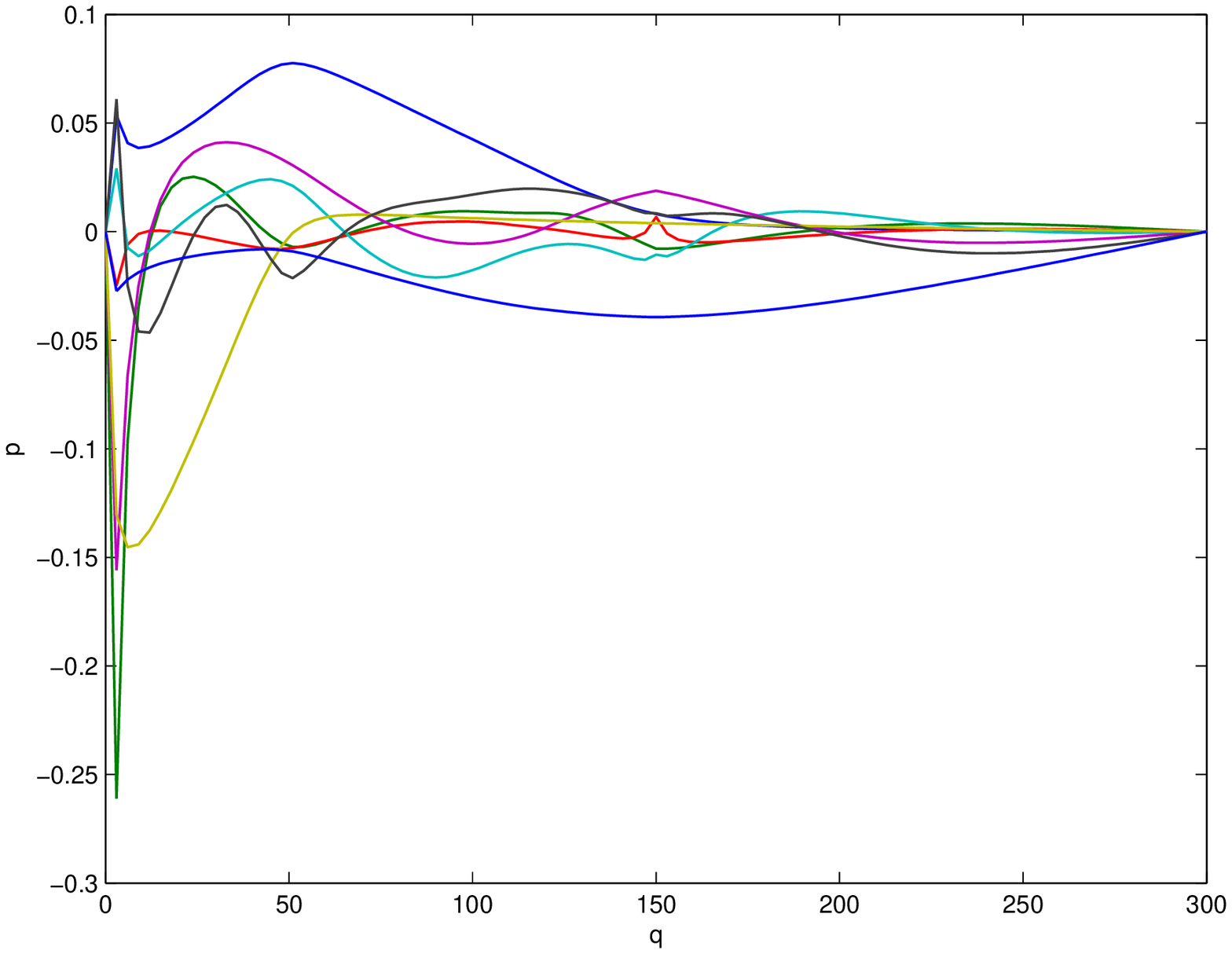}
\hfill
\psfrag{p}[c][t]{\footnotesize $\xi_j(s)$}
\includegraphics[width=.26\linewidth]{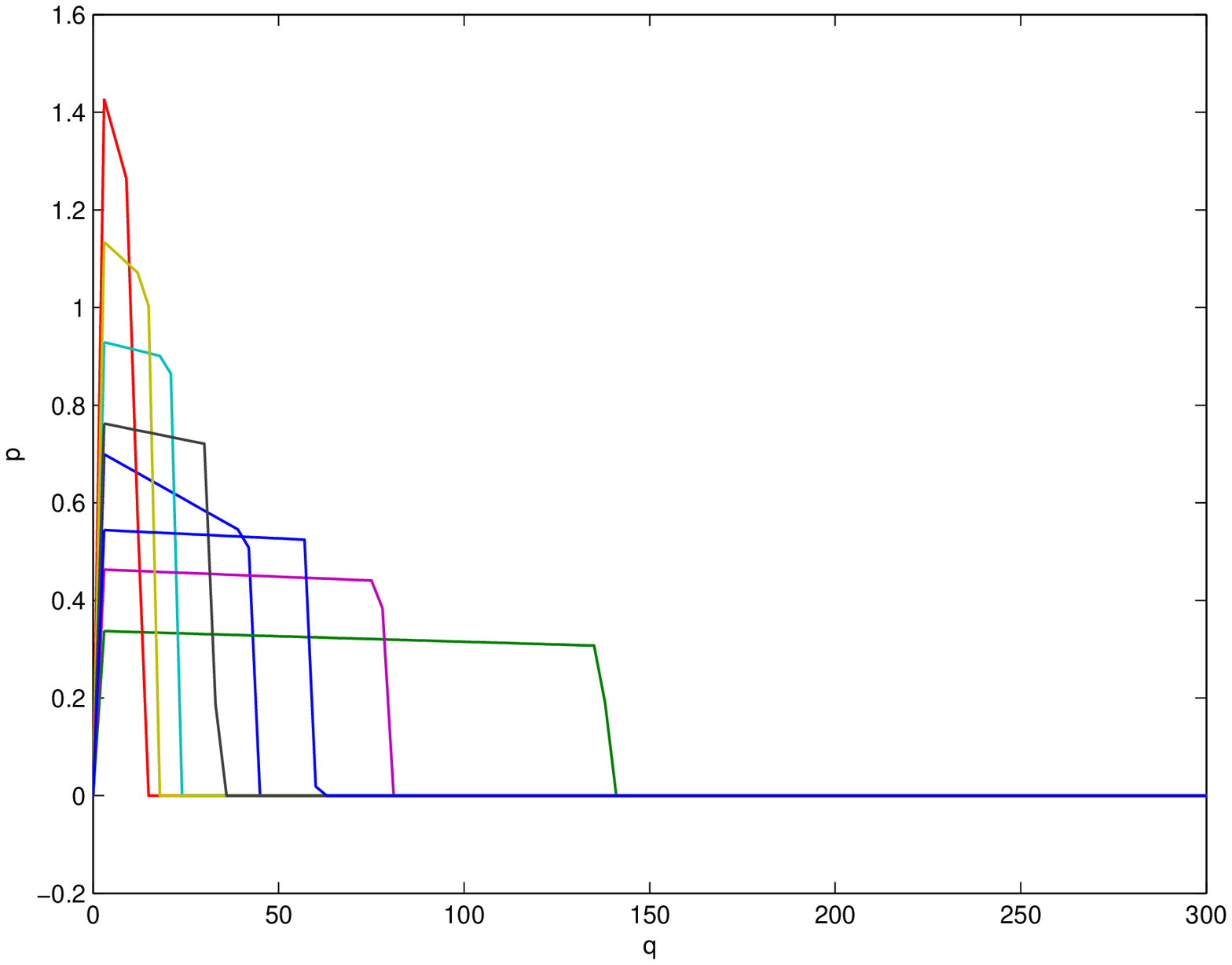}
\hfill
\psfrag{0}[r][b]{0}
\psfrag{p}[c][t]{\footnotesize $B\xi_j(s)$}
\includegraphics[width=.26\linewidth]{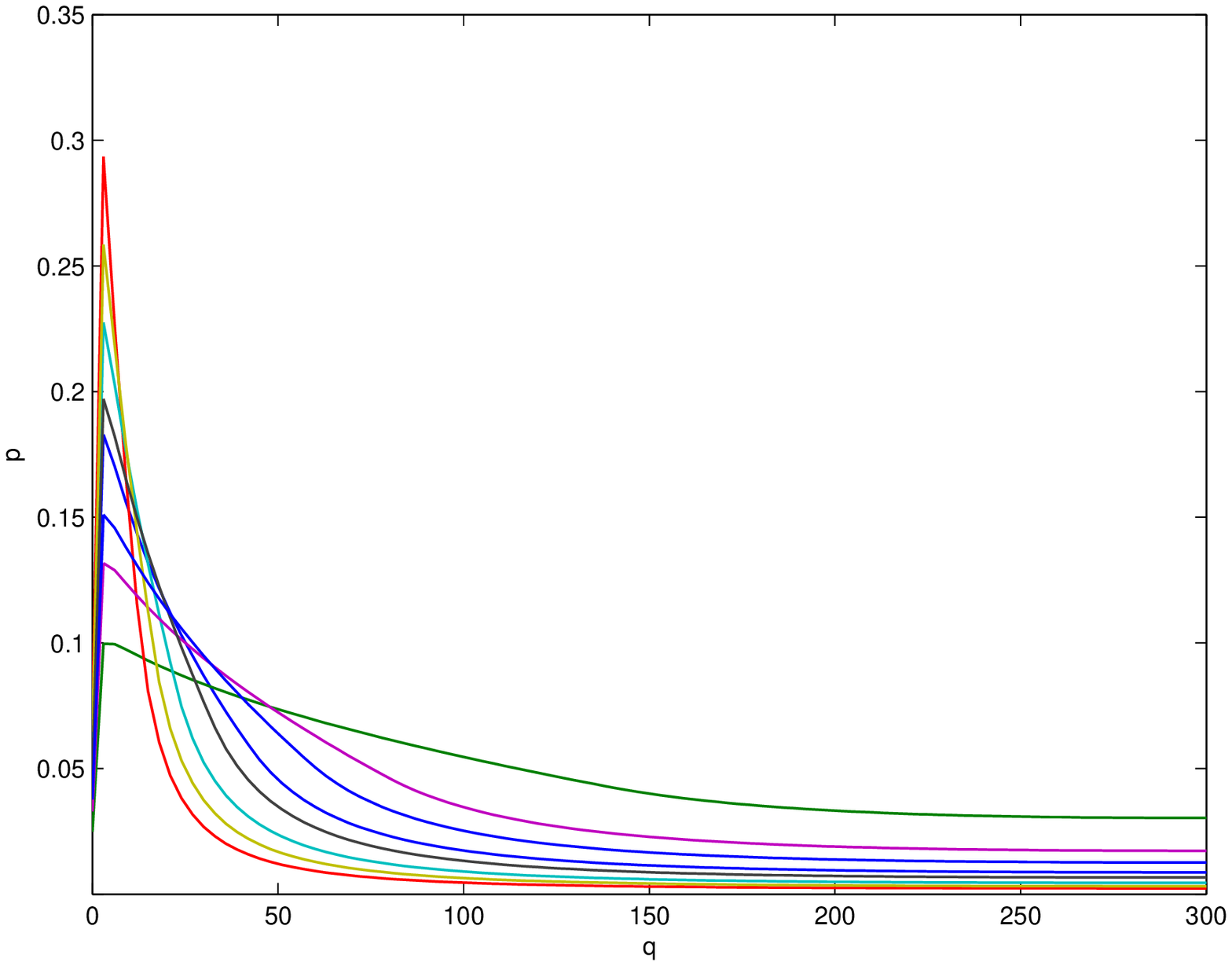}
\caption{Eight first vectors of the reduced basis $\Psi_N$, $\Xi_N$
  and the corresponding supremizers.}
\label{fig:basis}
\end{figure}

\begin{figure}[!h]
\centering
\psfrag{Time step=1}[c][t]{\footnotesize Time step =1}
\psfrag{Time step=10}[c][t]{\footnotesize Time step =10}
\psfrag{Time step=20}[c][t]{\footnotesize Time step =20}
\psfrag{0}[r][m]{\scriptsize 0}
\psfrag{150}[c][m]{\scriptsize 150}
\psfrag{300}[c][m]{\scriptsize 300}
\psfrag{40}[r][m]{\scriptsize 40}
\psfrag{80}[r][m]{\scriptsize 80}
\psfrag{20}[c][m]{}
\psfrag{60}[c][m]{}
\psfrag{50}[c][m]{}
\psfrag{100}[c][m]{}
\psfrag{200}[c][m]{}
\psfrag{250}[c][m]{}
\psfrag{p}[c][t]{\footnotesize $u_N$}
\psfrag{q}[c][b]{\footnotesize $s$}
\includegraphics[width=0.3\linewidth]{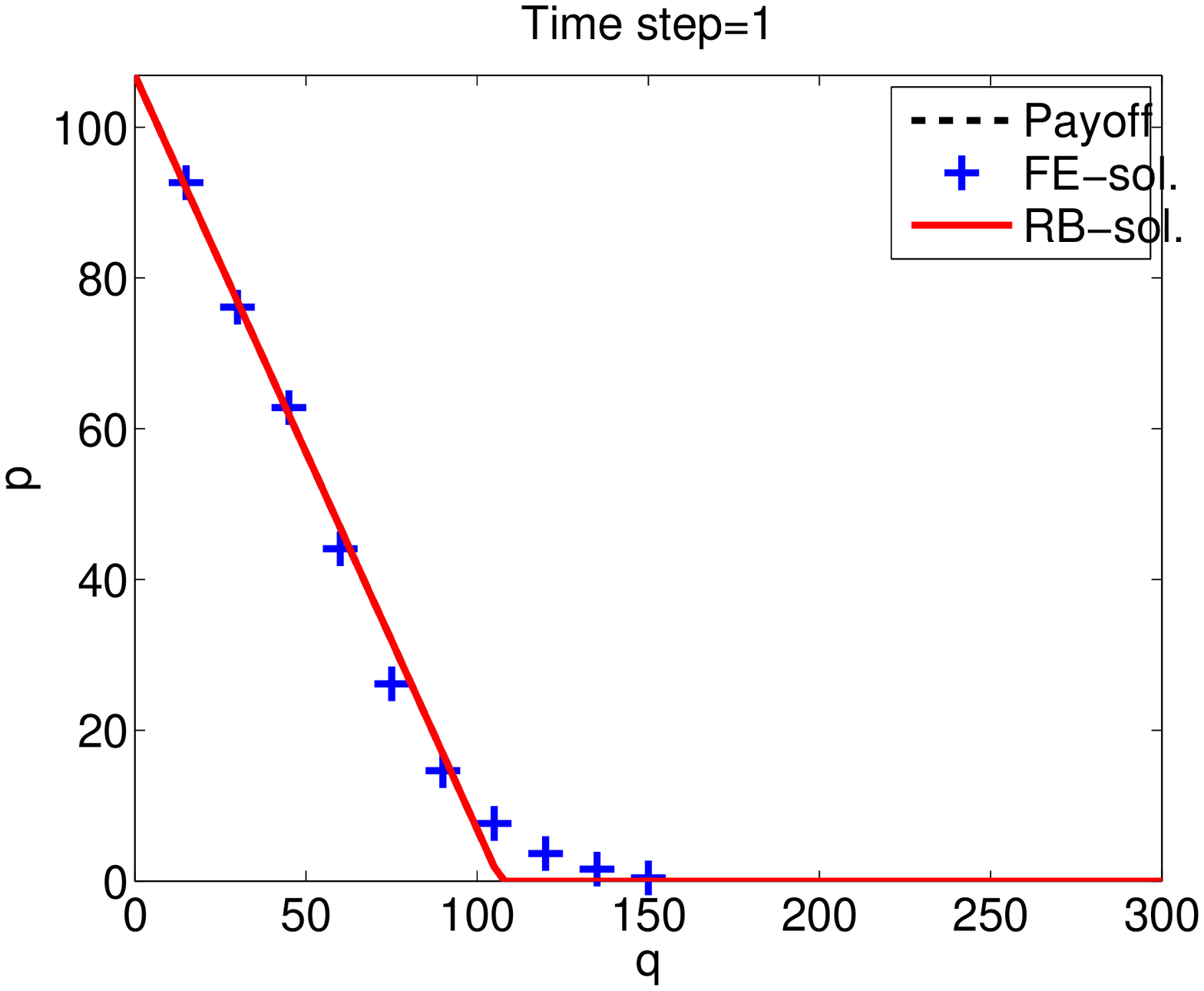}
\includegraphics[width=0.3\linewidth]{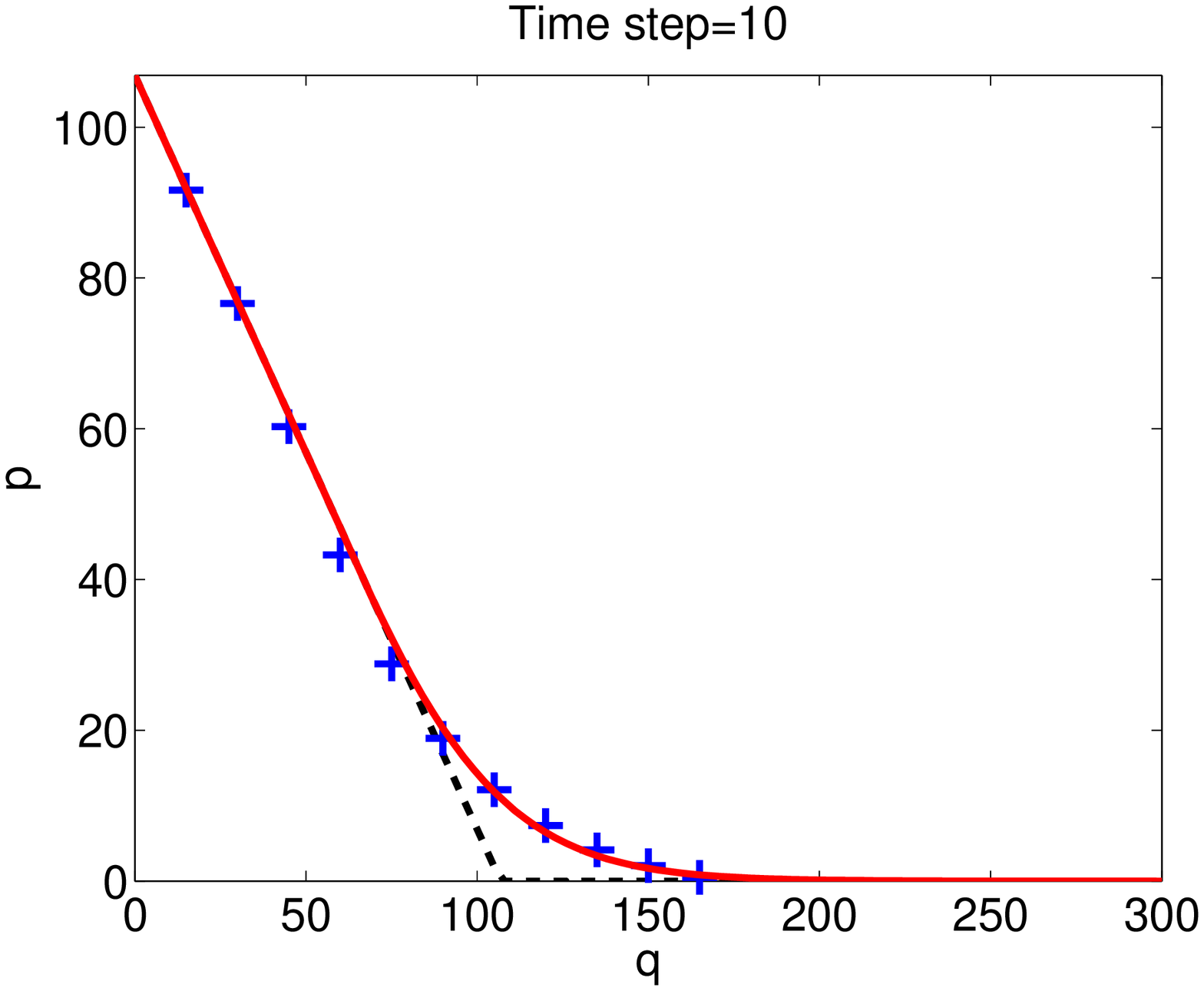}
\includegraphics[width=0.3\linewidth]{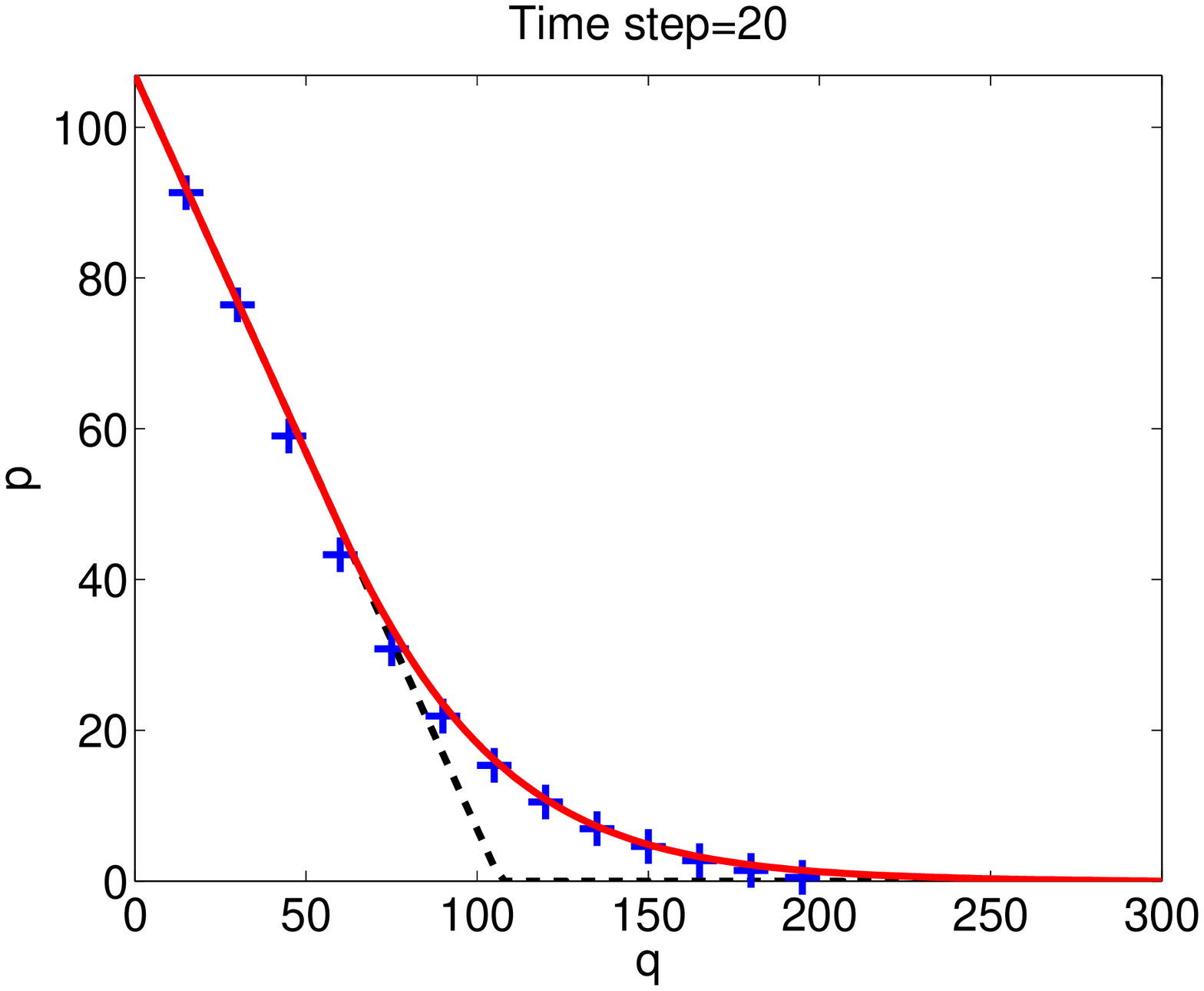}
\psfrag{Time step=1}[c][t]{}
\psfrag{Time step=10}[c][t]{}
\psfrag{Time step=20}[c][t]{}
\psfrag{0}[r][m]{\scriptsize 0}
\psfrag{150}[c][m]{\scriptsize 150}
\psfrag{300}[c][m]{\scriptsize 300}
\psfrag{40}[r][m]{\scriptsize 40}
\psfrag{80}[r][m]{\scriptsize 80}
\psfrag{20}[c][m]{}
\psfrag{60}[c][m]{}
\psfrag{50}[c][m]{}
\psfrag{100}[c][m]{}
\psfrag{200}[c][m]{}
\psfrag{250}[c][m]{}
\psfrag{p}[c][t]{\footnotesize $u_N$}
\psfrag{q}[c][b]{\footnotesize $s$}
\vspace{-.1cm}
\includegraphics[width=0.3\linewidth]{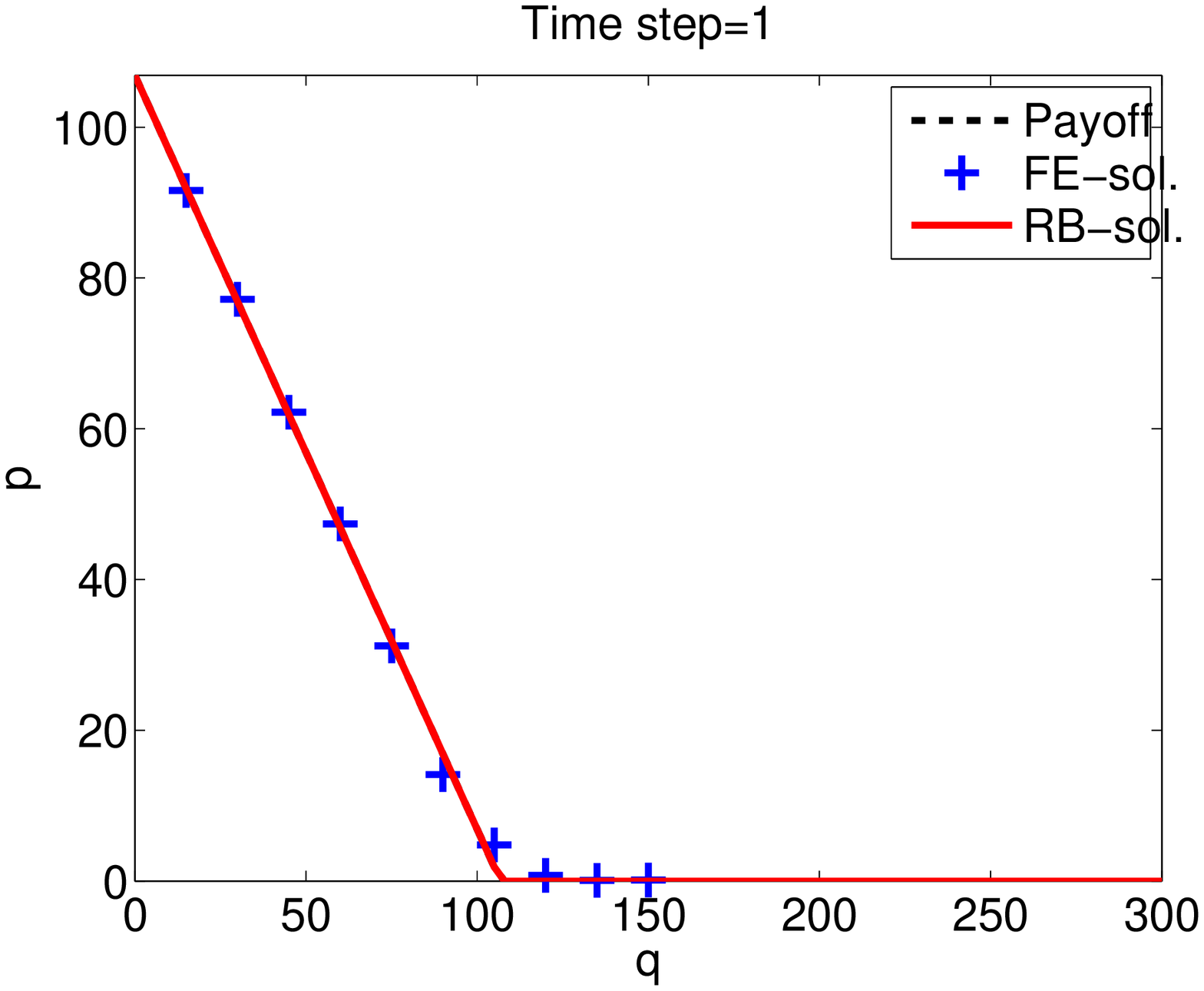}
\includegraphics[width=0.3\linewidth]{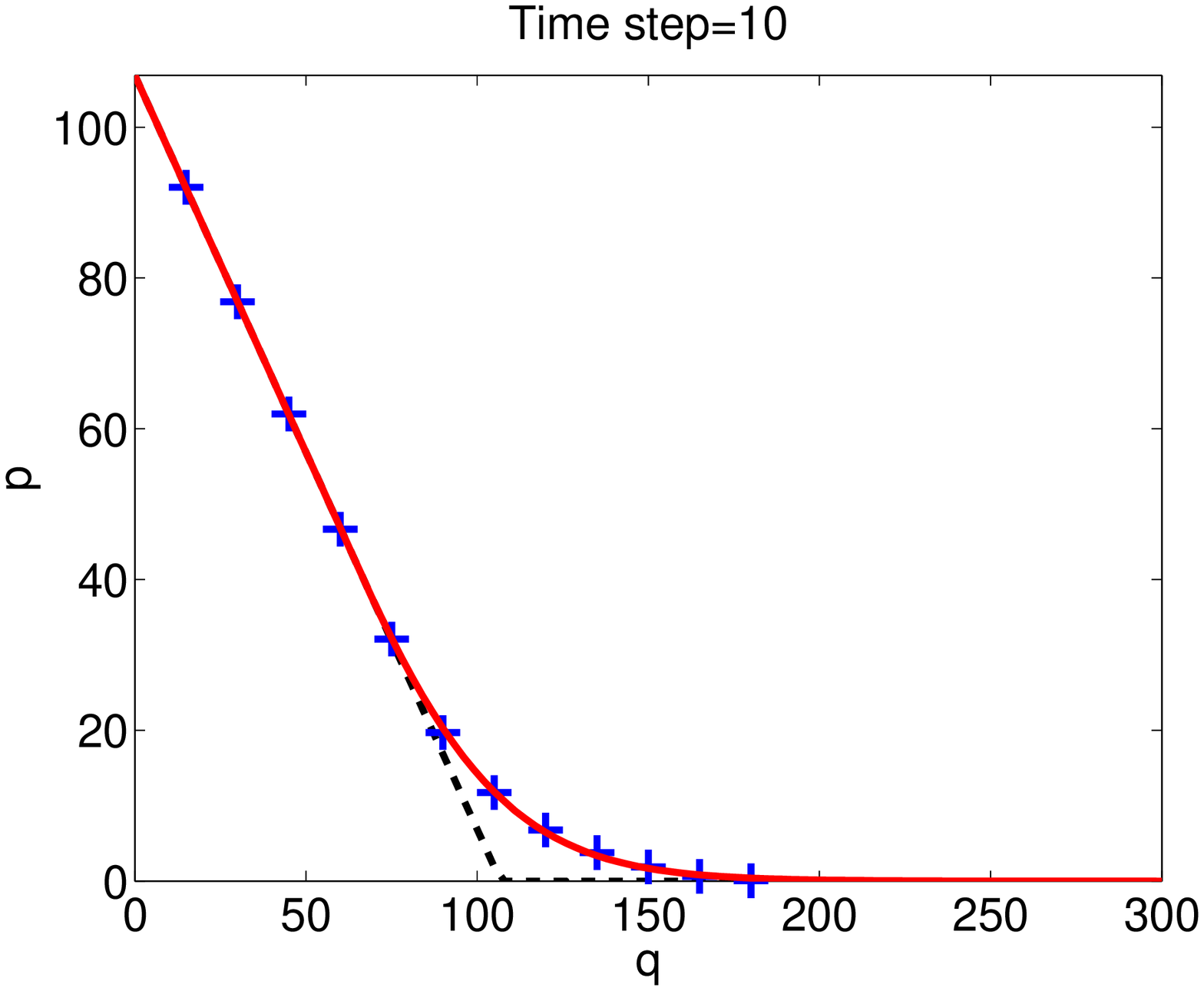}
\includegraphics[width=0.3\linewidth]{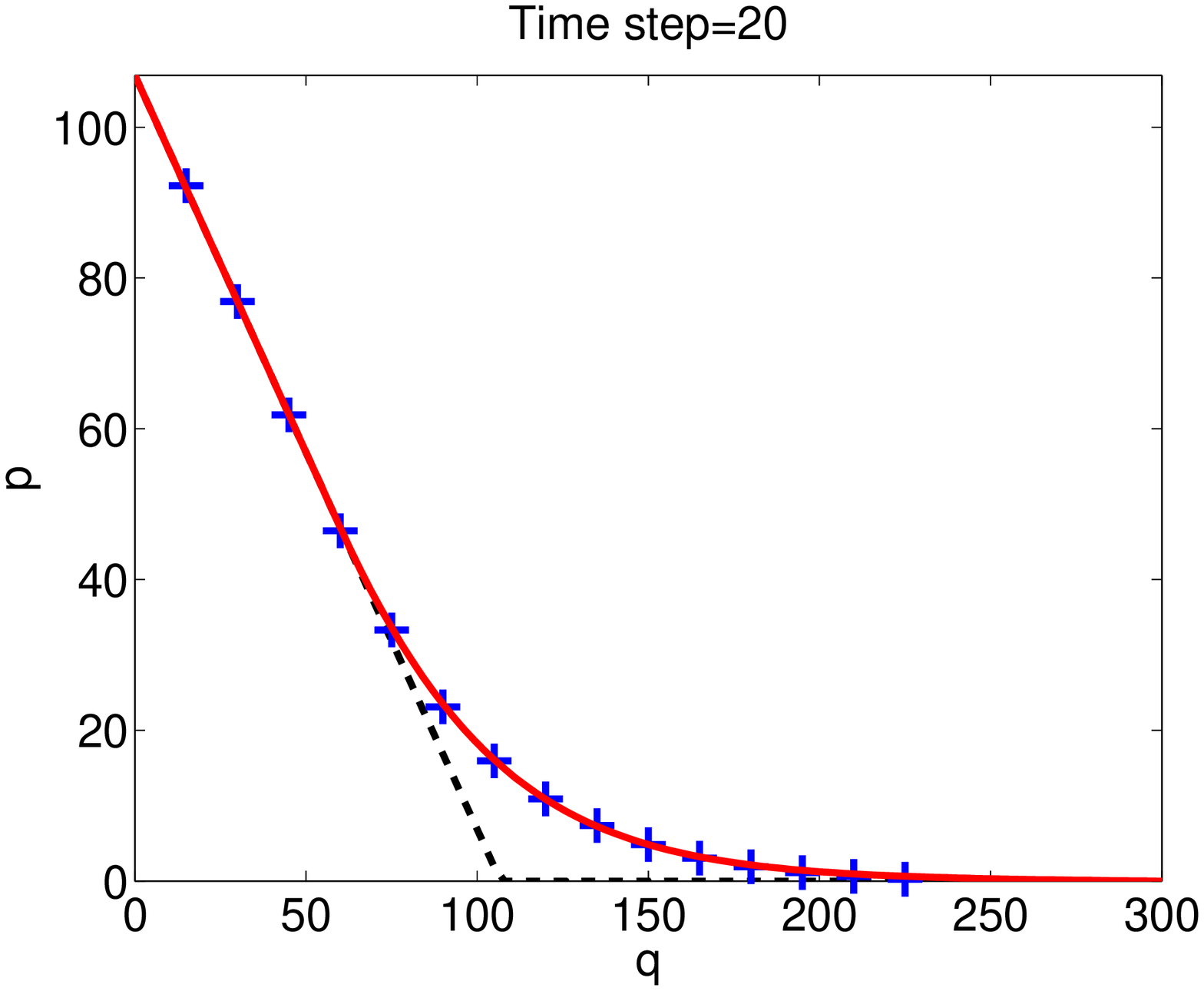}
\caption{Finite element approximation (solid red line) and Reduced basis
  approximation (blue $+$) at time steps $t/\Delta t=1$, $t/\Delta
  t=10$ and $t/\Delta t=T/\Delta t=20$. The payoff function $\psi$ is
  represented with the black dashed line. The reduced bases that are
  used have been generated by $(\widetilde N_V,N_W)=(8,8)$ (plots on the
  top) or $(\widetilde N_V,N_W)=(16,16)$ (plots on the
  bottom).}
\label{fig:simu}
\end{figure}

\begin{figure}
\centering
\psfrag{Angle-Greedy phi}[c][t]{\scriptsize Angle-Greedy}
\psfrag{POD-Greedy phi}[c][t]{\scriptsize POD-Greedy}
\psfrag{Max error}[c][t]{\scriptsize Max error}
\psfrag{10+3}[r][t]{\tiny $10^{3}$}
\psfrag{10+2}[r][t]{\tiny $10^{2}$}
\psfrag{10+1}[r][t]{\tiny $10^{1}$}
\psfrag{10-0}[r][t]{\tiny $1$}
\psfrag{10-1}[r][t]{\tiny $10^{-1}$}
\psfrag{10-3}[r][t]{\tiny $10^{-3}$}
\psfrag{10-5}[r][t]{\tiny $10^{-5}$}
\psfrag{10-7}[r][t]{\tiny $10^{-7}$}
\psfrag{1e 3}[c][t]{\tiny $10^3\phantom{11}$}
\psfrag{1e 2}[c][t]{\tiny $10^2\phantom{11}$}
\psfrag{1e 1}[c][t]{\tiny $10^1\phantom{11}$}
\psfrag{1e 0}[c][t]{\tiny $1$}
\psfrag{1e-1}[c][t]{\tiny $10^{-1}\phantom{111}$}
\psfrag{1e-2}[r][t]{\tiny $10^{-2}$}
\psfrag{1e-3}[r][t]{\tiny $10^{-3}$}

\psfrag{10-2}[r][t]{}
\psfrag{10-4}[r][t]{}
\psfrag{10-6}[r][t]{}
\psfrag{10-8}[r][t]{}
\psfrag{0}[c][b]{\tiny 0}
\psfrag{5}[c][b]{\tiny 5}
\psfrag{10}[c][b]{\tiny 10}
\psfrag{15}[c][b]{\tiny 15}
\psfrag{20}[c][b]{\tiny 20}
\psfrag{30}[c][b]{\tiny 30}
\psfrag{40}[c][b]{\tiny 40}
\psfrag{50}[c][b]{\tiny 50}
\psfrag{60}[c][b]{\tiny 60}
\psfrag{2}[r][t]{}
\psfrag{4}[r][t]{}
\psfrag{6}[r][t]{}
\psfrag{8}[r][t]{}
\psfrag{12}[r][t]{}
\psfrag{14}[r][t]{}
\psfrag{16}[r][t]{}
\psfrag{18}[r][t]{}

\psfrag{p}[c][t]{\scriptsize $\varepsilon_N^u$}
\psfrag{q}[c][b]{\scriptsize $N_v$}

\includegraphics[width=.3\linewidth]{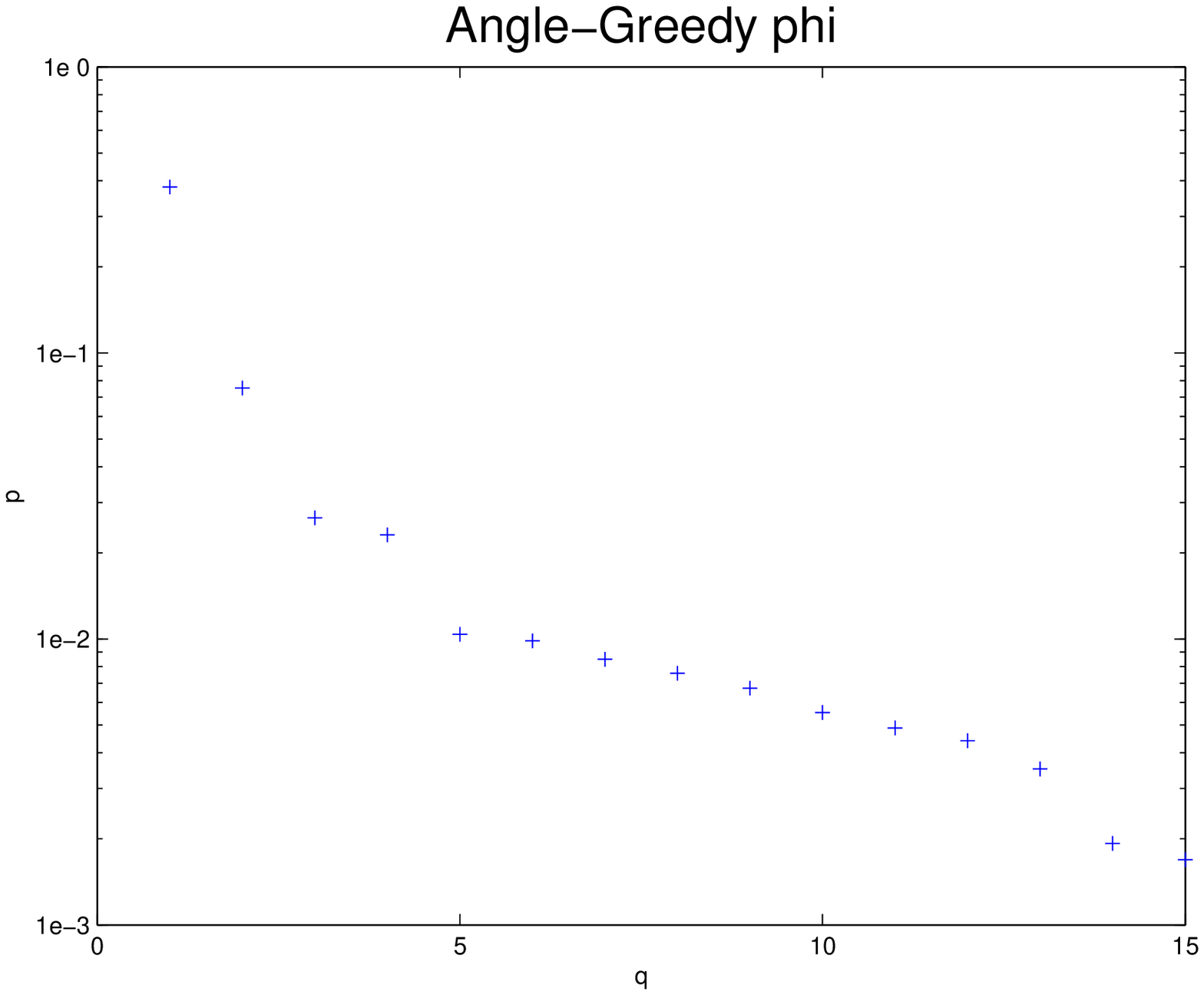}
\psfrag{p}[c][t]{\scriptsize $\begin{array}{c} \varepsilon_N^\lambda\\\phantom{1}\\\phantom{1}\end{array}$}
\psfrag{q}[r][b]{\scriptsize $N_W$}
\hfill
\includegraphics[width=.3\linewidth]{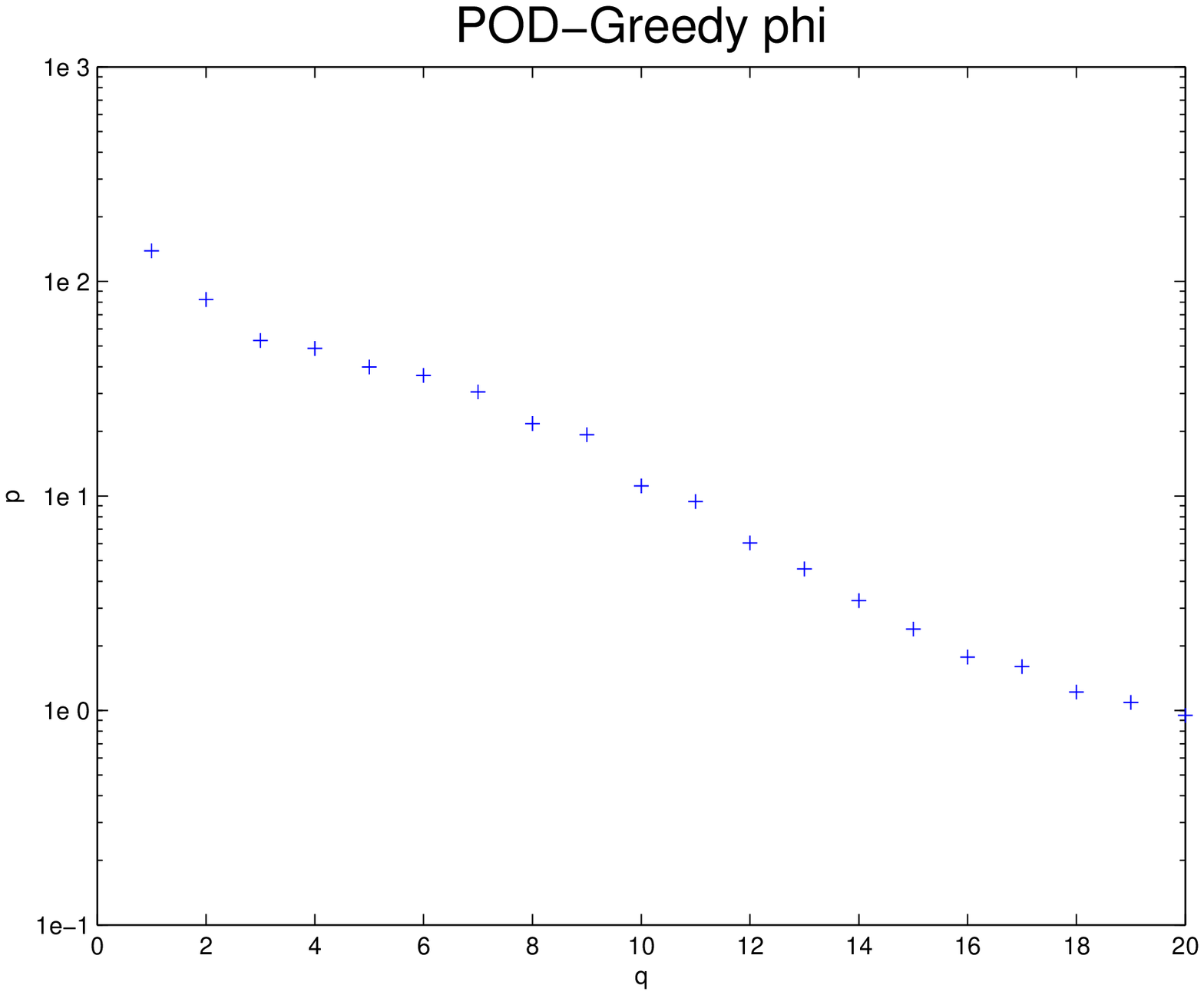}
\psfrag{a}[c][b]{\scriptsize $N_V$}
\psfrag{b}[c][b]{\scriptsize $N_W$}
\hfill
\psfrag{c}[c][t]{\scriptsize
$\begin{array}{c}Err^{L^\infty}_N\\\phantom{1}\end{array}$}
\includegraphics[width=.3\linewidth]{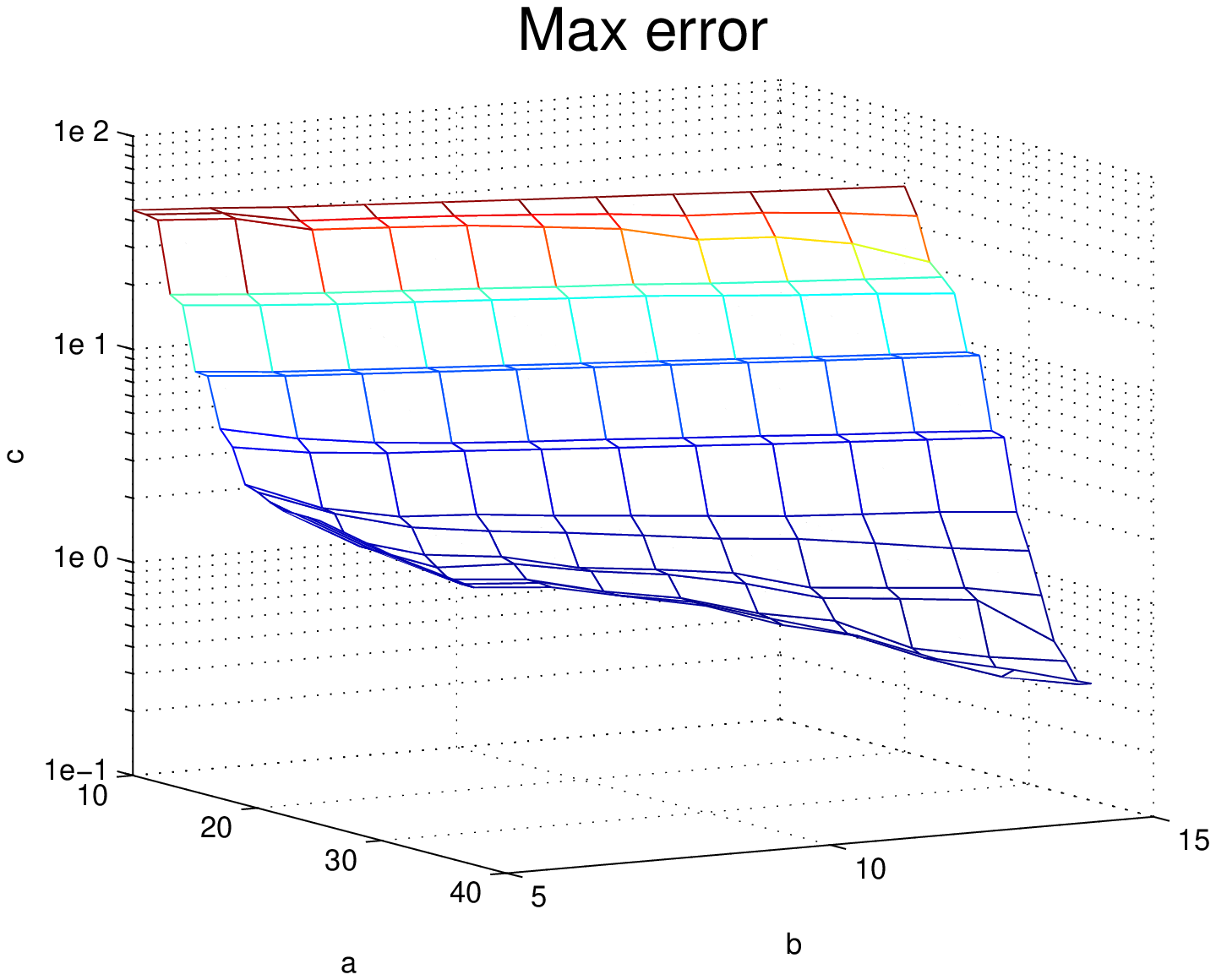}
\caption{Values of $\varepsilon_N^u$ and $\varepsilon_N^\lambda$ during the
  iterations of the greedy Algorithms~\ref{alg:ang-greedy} (left)
  and~\ref{alg:pod-greedy} (middle). Right: Values of   $Err^{L^\infty}_N$
with respect to $N_V$ and $N_W$.  }
\label{fig:ang-greedy}
\end{figure}

In this section, we present some numerical results obtained on the
American Option model. We start with a description of the numerical
values and methods we use.
%\subsection{Snapshots computation and basis construction}
\label{sec:RBdesign}
In order to compute snapshots, we use a standard finite element method
for the space discretization and the $\theta$-scheme presented in
Section~\ref{sec:2} for the
time-discretization. The time domain $[0,T]=[0,1]$ is discretized with a 
uniform mesh of step size $\Delta t:= T/L$, $L=20$. The
$\theta$-scheme is used with $\theta=1/2$, i.e. we apply a
Crank-Nicolson method.
The space domain $\Omega=(0,s_f)=(0,300)$ is discretized with a
uniform mesh of step size $\Delta s:=s_f/S$, $S=101$.  
For the function space, we use
standard conforming nodal first order finite elements. For the sake of simplicity, we keep the notation $V$ for the
discrete high dimensional space and define it by  
$V:= \{v \in H_0^1(\Omega)| v_{|[s_m,s_{m+1}]}\in
P_1,m=0,\ldots,S-1\}$ of dimension $H_V=H:=S-2=99$ with $s_m := m
\Delta s$. 
We associate the basis function $\phi_i \in V$ with its Lagrange 
node $s_i\in\Omega$, i.e., $\phi_i(s_j) = \delta_{ij},i,j=1,\ldots,H$.
The discretization of the Lagrange multipliers is performed using a dual finite
element basis $\chi_j$ of $W:=V'$ having the same support as
$\phi_j $, so that
$b(\phi_i,\chi_j)=\delta_{ij}$, $i,j=1,\ldots,H_W=H$. The cone $M$ is
defined by:
$M=\left\{\sum_{i=1}^{H_W}\eta_i\chi_i,\ \eta_i\geq 0 \right\}.$ To build the basis, we consider a subset ${\cal
  P}_{train}$ of $\cal P$  that is composed of $N=16$ values chosen randomly in the set 
\begin{eqnarray*}
{\cal  P}&=[(1-\frac \varepsilon 2)K_0,(1+\frac \varepsilon 2)K_0]
\times[(1-\frac \varepsilon 2)r_0,(1+\frac \varepsilon 2)r_0]\\&\times[(1-\frac \varepsilon 2)q_0,(1+\frac \varepsilon 2)q_0]
\times[(1-\frac \varepsilon 2)\sigma_0, (1+\frac \varepsilon 2)\sigma_0].
\end{eqnarray*}
with the numerical values $\varepsilon=0.1$, $K_0=100$,  $r_0=0.05$, $q_0=0.0015$, $\sigma_0=0.5$.
%\subsection{Example of trajectory simulation}
To define the
basis  $\Psi_N$ and the convex set  $\Xi_N$,  we use Algorithm~\ref{alg:pod-greedy}  combined with the
enlargement by the  supremizers and Algorithm~\ref{alg:ang-greedy}. The eight first vectors of $\Psi_N$, $\Xi_N$ and the supremizers are
represented in Figure~\ref{fig:basis}. We simulate two trajectories corresponding to the values
$(\widetilde N_V,N_W)=(8,8)$ and $(\widetilde N_V,N_W)=(16,16)$
respectively. The corresponding bases $\Psi_N$ are of size $N_V=16$
and $N_V=32$ respectively. We chose randomly a parameter 
vector $\mu$ corresponding to the 
values $K=106.882366$, $r=0.048470$, $d=0.007679$, $\sigma=0.418561$
in ${\cal  P}$. %A rough visual comparison between the finite element simulation and
%our reduced basis approximation does not enable to distinguish the two
%trajectories. 
Some steps of the simulation are represented in
Figure~\ref{fig:simu}, the top and lower row refer to 
the smaller and larger reduced spaces, respectively.
We clearly see the improvement in the approximation 
by increasing the reduced dimensions.
%\vspace{-1cm}
%\subsection{Test of Greedy algorithms}
In order to evaluate the efficiency of the greedy algorithms proposed
in Section~\ref{sec:RBcomp}, we plot the evolution of the
quantities 
\begin{eqnarray*}
\varepsilon_N^u := \max_{\mu\in{\cal P}_{train}}\sqrt{ \sum_{n=0}^L\|u^n(\mu)-\Pi_{V^k_N}(u^n(\mu))\|^2_V }, \
\varepsilon_N^\lambda := \max_{\footnotesize\two{n=0,\ldots,L ,} {\mu\in {\cal
    P}_{train}}} \left( \measuredangle \left( \lambda^n(\mu) ,W_N^k \right)\right)
\end{eqnarray*}
 during their iterations. 
The results are plotted in the first two diagrams in 
Figure~\ref{fig:ang-greedy}. We observe an excellent exponential convergence 
of the approximation measures.%\subsection{Efficiency of the method}
%In order to evaluate numerically the efficiency of our method,
As final experiment, we address the generalization ability of the RB-model to 
parameters outside the training set. 
We consider ${\cal P}_{test}\subset {\cal P}$, a random set of
$N_{test}=10$ parameter vectors and estimate, for a given $\mu\in{\cal P}$,  the efficiency of our method through these quantities:
$$ err_N(\mu)=\sqrt{\Delta t \sum_{n=0}^{L} \| u^n(\mu)-u^n_N(\mu)\|_V^2}, 
\quad Err_N^{L^\infty}=\max_{\mu \in {\cal P}_{test}}\left(err_N(\mu)\right).$$ 
Note that  $err_N(\mu)$ actually depends on $\Psi_N$ ; for the sake of
simplicity, we have omitted this reliance in the notation. As a test, we 
evaluate the influence of the parameters $\widetilde N_V,N_W$
determining the sizes of the bases $\Psi_N$ and $\Xi_N$ 
on $Err_N^{L^\infty}$.  The results are plotted
in the right diagram of Figure~\ref{fig:ang-greedy}. In our example 
we numerically obtain $N_V = \widetilde N_V + N_W$ in all cases, indicating, 
that the primal snapshots and supremizers are linearly independent.
We observe a reasonable good error decay when simultaneously
increasing $\widetilde N_V$ and $N_W$, indicating that the reduced method is
working well. We also note that in our case, the size of the dual
basis has a limited impact on the results. 
%\begin{figure} 
%\centering
%\psfrag{Angle-Greedy phi}[c][t]{\scriptsize Angle-Greedy}
%\psfrag{POD-Greedy phi}[c][t]{\scriptsize POD-Greedy}
%\psfrag{Max error}[c][t]{\scriptsize Max error}
%\psfrag{p}[c][t]{\footnotesize $\max_{\mu\in{\cal P}_{train}}{\left( \sum_{n=0}^L\|u^n(\mu)-\Pi_{V^k_N}(u^n(\mu))\|^2_V\right)}$}
%\psfrag{q}[c][t]{\footnotesize $N_v$}
%\includegraphics[width=.3\linewidth]{./outputEPS/decreasPODanglePrand_N=2}
%\psfrag{p}[c][t]{\footnotesize $\max_{\mu\in{\cal P}_{train}}{\left( \measuredangle \left( \lambda^n(\mu) , {\rm
%    span}\left\{\xi_1,\dots,\xi_{k}\right\} \right)\right)}$}
%\psfrag{q}[c][t]{\footnotesize $N_W$}
%\hfill
%\includegraphics[width=.3\linewidth]{./outputEPS/decreasPODgreedyPrand_N=2}
%\psfrag{a}[c][t]{\footnotesize $N_V$}
%\psfrag{b}[c][t]{\footnotesize $N_W$}
%\hfill
%\psfrag{c}[c][t]{\scriptsize $Err^{L^\infty}_N$}
%\includegraphics[width=.3\linewidth]{./outputEPS/RbError_NphiNximaxPrand_N=2}
%\caption{Values of $\varepsilon_N^u$ and $\varepsilon_N^\lambda$ during the 
%  iterations of the greedy Algorithms~\ref{alg:ang-greedy} (left)
%  and~\ref{alg:pod-greedy} (middle). Right: Values of   $Err^{L^\infty}_N$ with respect to $N_V$ and $N_W$.  }
%\label{fig:ang-greedy}
%\end{figure}

\vspace{-1cm}
\bibliographystyle{spbasic}
\bibliography{enumath}

%Moreover, we note some
%instabilities when $N_W$ is large than $N_V$ ; see the value of the
%error for $(N_V,N_W)=(10,15)$.

%\begin{acknowledgement}
%If you want to include acknowledgments of assistance and the like at
%the end of an individual chapter please use the \verb|acknowledgement|
%environment -- it will automatically render Springer's preferred
%layout. 
%\end{acknowledgement}
%

%\input{referenc}
\end{document}